\documentclass{article}

\usepackage{amsmath}
\usepackage{amssymb}
\usepackage{amsthm}
\usepackage{amscd}

\usepackage{amsxtra}
\usepackage{amsfonts}
\usepackage{eucal}
\usepackage{amssymb}

\newtheorem{theorem}{Theorem}[section]
\newtheorem{lemma}[theorem]{Lemma}
\newtheorem{conjecture}[theorem]{Conjecture}
\newtheorem{proposition}[theorem]{Proposition}

\theoremstyle{definition}

\newtheorem{example}[theorem]{Example}

\theoremstyle{remark}
\newtheorem{remark}[theorem]{Remark}

\numberwithin{equation}{section}

\makeatletter
\def\eqnarray{%
  \stepcounter{equation}%
  \let\@currentlabel=\theequation
  \global\@eqnswtrue
  \global\@eqcnt\z@
  \tabskip\@centering
  \let\\=\@eqncr
  $$\halign to \displaywidth\bgroup\@eqnsel\hskip\@centering
  $\displaystyle\tabskip\z@{##}$&\global\@eqcnt\@ne
  \hfil$\displaystyle{{}##{}}$\hfil
  &\global\@eqcnt\tw@$\displaystyle\tabskip\z@{##}$\hfil
  \tabskip\@centering&\llap{##}\tabskip\z@\cr}
\makeatother

\newfont{\germ}{eufm10}
\newfont{\slsmall}{cmsl8}

\def\La{\Lambda}
\def\la{\lambda}

\def\ol#1{\overline{#1}}

\def\veps{\varepsilon}

\def\C{{\mathbb C}}
\def\Z{{\mathbb Z}}

\begin{document}

\title{Difference $L$ operators related to $q$-characters}

\author{
A. Kuniba\thanks{Institute of Physics, University of Tokyo, Tokyo 153-8902, Japan},
\hskip0.2cm
M. Okado\thanks{Department of Informatics and Mathematical Science,
    Graduate School of Engineering Science,
Osaka University, Osaka 560-8531, Japan},
\hskip0.2cm
J. Suzuki\thanks{Department of Physics,
Faculty of Science,
Shizuoka University, Ohya 836, Japan},
\hskip0.2cm
Y. Yamada\thanks{Department of Mathematics, Faculty of Science,
Kobe University, Hyogo 657-8501, Japan}
}
\date{}
\pagestyle{plain}

\maketitle
\begin{abstract}
We introduce a factorized difference operator $L(u)$ 
annihilated by the 
Frenkel-Reshetikhin screening operator for 
the quantum affine algebra $U_q(C^{(1)}_n)$.
We identify the coefficients of $L(u)$ with the 
fundamental $q$-characters, and establish a number of formulas
for their higher analogues.
They include Jacobi-Trudi and 
Weyl type formulas, canceling tableau sums,  
Casorati determinant solution to the $T$-system, and so forth.
Analogous operators for the orthogonal series 
$U_q(B^{(1)}_n)$ and $U_q(D^{(1)}_n)$ are also presented.
\end{abstract}

\section{Introduction}\label{sec:intro}

In this paper we introduce a factorized difference operator 
$L(u)$ related to the quantum affine algebra $U_q(C^{(1)}_n)$,
and present its application to the $q$-characters 
of some finite dimensional representations.
As for basic facts on finite dimensional representations, 
we refer to \cite{CP1,CP2} and \cite{Ka}.

The theory of $q$-characters was introduced in 
\cite{FR2,FM} motivated by their study of 
deformed ${\mathcal W}$-algebras.
See also \cite{Kn}.
The $q$-characters $\chi_q$ are 
Laurent polynomials in infinitely many variables 
$\{Y_a(u)^{\pm 1} \mid 1 \le a \le n, u \in \C\}$, 
which reduce to linear combinations of 
the usual characters with respect to a classical 
simple subalgebra in the limit $q \rightarrow 1$.
For an irreducible representation $V$,
$\chi_q(V)$  contains the highest weight monomial with coefficient 1, and 
the rest is generated by multiplying lowering factors 
corresponding to the negative roots.
One of the fundamental properties of the 
$q$-characters is that they enjoy a symmetry analogous to 
simple reflections in the Weyl group  \cite{FR2,FM}.
It is represented as 
$S_a \cdot \chi_q = 0$ for all $1 \le a \le n$, 
where $S_a$ is called the screening operator.

We construct a difference operator 
$L(u)$ acting on functions of a variable $u$ 
with the following features:

\vskip0.2cm
\noindent
$\bullet$  $S_a \cdot L(u) = 0$ for all $a$.

\noindent
$\bullet$ $L(u)$ generates all the 
fundamental $q$-characters of $U_q(C^{(1)}_n)$.

\vskip0.2cm
By fundamental $q$-characters we mean those 
for the fundamental representations in the sense of \cite{CP1}.
Such an  operator of order $n$ was known for 
$A^{(1)}_n$ in \cite{FR1,FRS}. 
Curiously, our $L(u)$ for $C^{(1)}_n$ is of order $2n+2$  and 
has a form similar to the $A^{(1)}_{2n+1}$ case, which contrasts 
with the well known embedding $C_n \hookrightarrow A_{2n-1}$.
It contains each fundamental $q$-character twice and 
admits a factorization into $2n$ product of first order 
difference operators and a second order one.
The last piece is further factorized if one 
extends the coefficient field 
${\mathcal Y} = \Z[Y_a(u)^{\pm 1}]_{1 \le a \le n, u \in \C}$ to 
$\Z[Q_a(u)^{\pm 1}]_{1 \le a \le n, u \in \C}$ by 
introducing the Baxter $Q$ functions as (\ref{eq:defy}).
Such an identification 
is based on the connection \cite{FR2} between 
$q$-characters and the analytic Bethe ansatz \cite{R,KS,KOS} 
for solvable lattice models \cite{B}.

To utilize the  ideas in the latter is another 
motivation of the paper.
Roughly, $\chi_q$ corresponds to an eigenvalue 
formula for transfer matrices, and 
the condition $\chi_q \in \cap_a \mbox{Ker} S_a$ 
to its  pole-freeness.
Our approach here is based on the difference equation $L(u)w(u) =0$, 
and involves only elementary linear algebra and a portion 
of combinatorics.
By means of a difference analogue of the Wronskian method for 
ordinary linear differential equations, we express  
the fundamental $q$-characters and their higher analogues 
in terms of a ratio of Casorati determinants.
They may be viewed as analogues of the Weyl formula 
for usual characters.
By a standard argument \cite{NNSY}, 
the ratio of the Casorati determinants is 
equal to a sum over semistandard tableaux on letters 
$\{1, 2, \ldots, 2n+2\}$.
Until this point, parallel results are derivable also for 
the $A^{(1)}_{2n+1}$ case.
A curiosity for $C^{(1)}_n$ case is that one of the tableau variable
has a minus sign. (See (\ref{eq:xdef2}).)
Nevertheless all the negative contributions 
cancel out for many examples including fundamental $q$-characters as 
exhibited in Appendix \ref{app:proof}.
It will be interesting to seek applications of the present 
results in the light of the works
\cite{BHK,CK,DDT,FRS,KLWZ,M2,SS,S1}.

The paper is organized as follows.
In Section \ref{sec:L}, we define $L(u)$ and derive 
various formulas for the fundamental $q$-characters.
In Section \ref{sec:higher} we extend the results in 
the preceding section slightly to  
a class generated by $L(u)$.
This is partly motivated by a similar structure 
observed for the Stokes multipliers \cite{DDT,S1,S2}.
In Section \ref{sec:Tsys}, we present a solution of the 
$C^{(1)}_n$ $T$-system \cite{KNS} in terms of 
the Casorati determinants.
An analogous result for $A^{(1)}_n$ is available in \cite{KLWZ}.
In Section \ref{sec:BD} we 
give  similar difference operators $L(u)$ for 
$B^{(1)}_n$ and $D^{(1)}_n$, 
which originate in \cite{KOS,TK}.
However, their orders are not finite as opposed to the $C^{(1)}_n$ case.
Appendix \ref{app:proof} contains a combinatorial proof 
of Proposition \ref{pr:main}.
Appendix \ref{app:iran} provides some basic lemmas 
connecting the Casorati determinants, tableau sums 
and Jacobi-Trudi type formulas.

\section{\mathversion{bold} 
$L(u)$ and fundamental $q$-characters}\label{sec:L}

Let $\{\alpha_a\mid 1 \le a \le n\}$ and $\{\La_a \mid 1 \le a \le n\}$ 
be the sets of simple roots and fundamental weights of 
$C_n$ normalized as $(\alpha_a \vert \alpha_a) = 1+\delta_{a n}$.
Throughout the paper we set 
\[ N = 2n+2. \]
We denote the Frenkel-Reshetikhin variable 
$Y_{a,q^u}$ \cite{FR2} for $U_q(C^{(1)}_n)$
by $Y_a(u)\; (1 \le a \le n)$, and introduce the
Baxter $Q$ functions $Q_a(u)$ related to $Y_a(u)$ as 
\begin{equation}\label{eq:defy}
Y_a(u) = \frac{Q_a(u-\frac{(\alpha_a\vert \alpha_a)}{2})}
{Q_a(u+\frac{(\alpha_a\vert \alpha_a)}{2})}.
\end{equation}
For the definition of the $q$-character we refer to 
the original paper \cite{FR2}.
Here we recall the screening operator $S_a$.
$S_a$ sends the elements in 
${\mathcal Y} = \Z[Y_a(u)^{\pm1}]_{1 \le a \le n, u \in \C}$ to the 
ring extended by adjoining the extra symbols 
$\{S_a(u) \mid 1 \le a \le n, u \in \C\}$.
The action is given by 
\[
S_a \cdot Y_b(u) = \delta_{a b} Y_b(u)S_b(u)
\] 
and the Leibniz rule. 
Thus for example 
${S}_a \cdot Y_b(u)^{-1} = -\delta_{a b} Y_b(u)^{-1}S_b(u)$.
The symbol $S_a(u)$ is  assumed to obey the following 
relation in the extended ring:
\begin{equation}\label{eq:SAS}
\begin{split}
&S_a\bigl(u+(\alpha_a\vert \alpha_a)\bigr) 
= A_a\left(u+\frac{(\alpha_a\vert \alpha_a)}{2}\right)S_a(u),\\
&A_a(u) = \prod_{b=1}^{n}
\frac{Q_b(u-(\alpha_a\vert\alpha_b))}{Q_b(u+(\alpha_a\vert\alpha_b))},
\end{split}
\end{equation}
where $A_a(u)$ can actually be expressed in terms of the $Y$-variables only.

Let 
\[  
J = \{1 \prec 2 \prec \cdots \prec n \prec \bar{n} \prec
\cdots \prec \bar{2} \prec \bar{1} \}
\]
be an ordered set.
In what follows we shall work with the 
two sets of variables 
$\{z_a(u) \mid a \in J\}$ and $\{ x_a(u) \mid 1 \le a \le N\}$ 
specified by
\begin{align}
&z_a(u) = \frac{Y_a(u+\frac{a}{2})}{Y_{a-1}(u+\frac{a+1}{2})}, \quad 
z_{\bar{a}}(u) = \frac{Y_{a-1}(u+\frac{2n-a+3}{2})}{Y_{a}(u+\frac{2n-a+4}{2})}
\quad 1 \le a \le n, \label{eq:defz}\\
&x_a(u) = z_a(u),\quad x_{2n+3-a}(u) = z_{\bar a}(u) \quad 1 \le a \le n,
\label{eq:xdef1}\\
&x_{n+1}(u) = - x_{n+2}(u) = 
\frac{Q_n(u+\frac{n}{2})Q_n(u+\frac{n+4}{2})}{Q_n(u+\frac{n+2}{2})^2},
\label{eq:xdef2}
\end{align}
where $Y_0(u) = 1$.
Note here that unlike (\ref{eq:defz}) and (\ref{eq:xdef1}), 
$x_{n+1}(u)= -x_{n+2}(u)$ is {\em not } expressible 
as  a ratio of $Y_n$'s.

In this paper, we denote the $q$-character $\chi_q(V_{\La_a}(q^u))$ 
of the fundamental representation $V_{\La_a}(q^u)$ by
$T^{(a)}_1(u)\; (1 \le a \le n)$.
In terms of the variables (\ref{eq:defz}), 
it is given by
\begin{equation}\label{eq:Tdef1}
T^{(a)}_1(u+\frac{1}{2}) = \sum z_{i_1}(u+\frac{a-1}{2})
z_{i_2}(u+\frac{a-3}{2}) \cdots 
z_{i_a}(u-\frac{a-1}{2}),
\end{equation}
where the sum runs over $i_1, \ldots, i_a \in J$ such that 
\begin{align}
&i_1 \prec \cdots \prec i_a,\label{eq:cond1}\\
&\mbox{ if } i_k = c,\; i_l= \bar{c}\, \mbox{ for some } 1 \le c \le n, 
\mbox{ then } n+k-l \ge c. \label{eq:cond2}
\end{align}
Upon a suitable convention adjustment, 
this agrees with $s_a(zq^{-1/2})$ in section 11.3 of \cite{FR1}.
It also coincides with 
$\La^{(a)}_1(u+1/2)$ in \cite{KS}, where 
the condition (\ref{eq:cond1})--(\ref{eq:cond2}) was 
first introduced under which 
the number of summands is indeed equal to the dimension 
$\binom{2n}{a} - \binom{2n}{a-2}$ of the $a$ th
fundamental representation of $C_n$.
We have $T^{(a)}_1(u) = Y_a(u)+ \cdots$, where $Y_a(u)$ is the highest
weight monomial \cite{FM}. 
\begin{example}
For $n=2$, the two fundamental $q$-characters read
\begin{align*}
T^{(1)}_1(u) &=
Y_1(u) + \frac{Y_2(u+\frac{1}{2})}{Y_1(u+1)} +
\frac{Y_1(u+2)}{Y_2(u+\frac{5}{2})}+ \frac{1}{Y_1(u+3)},\\
T^{(2)}_1(u) &=Y_2(u) + 
\frac{Y_1(u+\frac{1}{2})Y_1(u+\frac{3}{2})}{Y_2(u+2)}+
\frac{Y_1(u+\frac{1}{2})}{Y_1(u+\frac{5}{2})}+
\frac{Y_2(u+1)}{Y_1(u+\frac{3}{2})Y_1(u+\frac{5}{2})}+
\frac{1}{Y_2(u+3)}.
\end{align*}
\end{example}

Let $D$ be a difference operator $Dg(u) = g(u+1)D$.
We use the notation 
$\displaystyle\stackrel{\longrightarrow}{\prod_{i=1}^k}X_i=X_1X_2\cdots X_k$, 
$\displaystyle\stackrel{\longleftarrow}{\prod_{i=1}^k}X_i=X_k\cdots X_2X_1$.
By a direct calculation one finds
\begin{lemma}\label{lem:kakikae}
{\small 
\begin{align*}
&\displaystyle\stackrel{\longrightarrow}{\prod_{a=1}^n}
(1-z_{\bar{a}}(u)D)\cdot (1-z_{\bar{n}}(u)z_n(u+1)D^2)\cdot
\displaystyle\stackrel{\longleftarrow}{\prod_{a=1}^n}(1-z_a(u)D)\\
&= \displaystyle\stackrel{\longrightarrow}{\prod_{a=1}^n}
(z_{a}(u+n+1-a)-D)\cdot (z_{\bar{n}}(u-1)z_n(u)-D^2)\cdot
\displaystyle\stackrel{\longleftarrow}{\prod_{a=1}^n}(z_{\bar{a}}(u-n-2+a)-D)
\\
&= \displaystyle\stackrel{\longleftarrow}{\prod_{i=1}^N}
(1-\epsilon_ix_i(u)D)
= - \displaystyle\stackrel{\longrightarrow}{\prod_{i=1}^N}
(\epsilon_ix_i(u+n+1-i)-D),
\end{align*}}
where $\epsilon_i=1$ for $i \neq n+1,n+2$ and 
$\epsilon_{n+1}=\epsilon_{n+2} = \pm 1$.
\end{lemma}
In particular the middle quadratic factor can be factorized as
\[
D^2-z_{\bar{n}}(u-1)z_n(u)
= D^2-x_{n+1}(u-1)x_{n+1}(u) = (x_{n+1}(u) \pm D)(x_{n+2}(u-1) \pm D).
\]
We define the $L$ operator by
\begin{equation}\label{eq:Ldef}
L(u) = \displaystyle\stackrel{\longrightarrow}{\prod_{i=1}^N}
(x_i(u+n+1-i)-D) 
= - \displaystyle\stackrel{\longleftarrow}{\prod_{i=1}^N}
(1-x_i(u)D).
\end{equation}
Due to Lemma \ref{lem:kakikae}, $L(u)$ 
is a polynomial in $D$ of order $N$ with 
coefficients in ${\mathcal Y}$.
For a difference operator of the form 
$\sum_j c_j(u) D^j$ with $c_j(u) \in {\mathcal Y}$, 
the screening operator acts as
\[
S_a \cdot 
\Bigl(\sum_j c_j(u) D^j\Bigr)
= \sum_j \bigl(S_a\cdot c_j(u) \bigr)D^j.
\]

\begin{proposition}\label{pr:Lproperty}
\begin{align}
& L(u)Q_1(u) = 0, \label{eq:LQ}\\
&{S}_a \cdot L(u) = 0 \quad 1 \le a \le n. \label{eq:SL}
\end{align}
\end{proposition}
\begin{proof}
The rightmost factor in the first expression of (\ref{eq:Ldef}) reads
$(Q_1(u+1)/Q_1(u)-D)$, proving (\ref{eq:LQ}).
As for (\ref{eq:SL}), we illustrate the $a=n$ case.
By the definition, 
$S_n$ acts non-trivially only on the middle three factors in 
(\ref{eq:Ldef}), which is expanded as $(v = u+\frac{n}{2})$
\begin{equation*}
\begin{split}
&\frac{Y_{n-1}(v-\frac{1}{2})}{Y_{n-1}(v+\frac{3}{2})}-
\left(\frac{Y_{n-1}(v+\frac{1}{2})}{Y_{n}(v+2)} + 
\frac{Y_{n}(v)}{Y_{n-1}(v+\frac{3}{2})}\right)D \\
&+\left(\frac{Y_{n-1}(v+\frac{5}{2})}{Y_{n}(v+3)} + 
\frac{Y_{n}(v+1)}{Y_{n-1}(v+\frac{3}{2})}\right)D^3 + D^4.
\end{split}
\end{equation*}
Upon applying $S_a$, this becomes
\begin{equation*}
\begin{split}
&\left(S_n(v+2)\frac{Y_{n-1}(v+\frac{1}{2})}{Y_{n}(v+2)} -
S_n(v)\frac{Y_{n}(v)}{Y_{n-1}(v+\frac{3}{2})}\right)D \\
&+\left(-S_n(v+3)\frac{Y_{n-1}(v+\frac{5}{2})}{Y_{n}(v+3)} + 
S_n(v+1)\frac{Y_{n}(v+1)}{Y_{n-1}(v+\frac{3}{2})}\right)D^3.
\end{split}
\end{equation*}
This vanishes under (\ref{eq:SAS}), i.e.,
\[
S_n(v+2) = \frac{Y_n(v)Y_n(v+2)}
{Y_{n-1}(v+\frac{3}{2})Y_{n-1}(v+\frac{1}{2})}S_n(v).
\]
\end{proof}

Let us introduce the notation:
\begin{align}
&X_u(i_1,\ldots, i_a) = \prod_{k=1}^ax_{i_k}(u+1-k),\quad 
Z_u(i_1,\ldots, i_a) = \prod_{k=1}^az_{i_k}(u+1-k),\nonumber\\
&\sum_i X_u(i_1,\ldots, i_a) = \sum_{1 \le i_1 < \cdots < i_a \le N}
X_u(i_1,\ldots, i_a),\nonumber \\
&\sum_i Z_u(i_1,\ldots, i_a) = 
\sum_{1 \preceq i_1 \prec \cdots \prec i_a \preceq \bar{1}}
Z_u(i_1,\ldots, i_a).\label{eq:zsum}
\end{align}
Expanding  (\ref{eq:Ldef})  we get  the two expressions:
\begin{equation*}\label{eq:2exp}
\begin{split}
L(u) &= \sum_{a=0}^{N}(-1)^a \left(
\sum_i X_{u+n}(i_1,\ldots,i_{N-a}) \right) D^a,\\
&= -\sum_{a=0}^{N}(-1)^a \left(
\sum_i X_{u+a-1}(i_1,\ldots,i_{a}) \right) D^a.
\end{split}
\end{equation*}
Thus we find 
\begin{align}
L(u) &= \sum_{i=n+2}^{N}(-1)^ie_{N-i}(u+\frac{i}{2})D^i - 
\sum_{i=0}^n(-1)^ie_i(u+\frac{i}{2})D^i,\label{eq:Le}\\
&e_a(u+\frac{a}{2}) = -e_{N-a}(u+\frac{a}{2}) \quad 0 \le a \le N,
\label{eq:esym}\\
&e_a(u+\frac{a}{2}) := \sum_i X_{u+a-1}(i_1,\ldots,i_{a}).\label{eq:edef}
\end{align}
In particular (\ref{eq:esym}) tells $e_{n+1}(u) = 0, 
e_{N}(u) = -1$.
In view of (\ref{eq:xdef2}), the sum (\ref{eq:edef}) contains 
sign factors. 
However, they all cancel out leaving the 
``positive" contributions only.

\begin{proposition}[Canceling tableau sums for fundamental $q$-characters]
\label{pr:main}
We have $T^{(a)}_1(u) = e_a(u)$, namely, 
\begin{equation}\label{eq:tx}
T^{(a)}_1(u) = 
\sum_i X_{u+\frac{a}{2}-1}(i_1,\ldots,i_{a})\quad 1 \le a \le n.
\end{equation}
\end{proposition}
The proof is available in Appendix \ref{app:proof}, where we show that the
cancellation in (\ref{eq:edef}) precisely 
leaves the sum in (\ref{eq:Tdef1})--(\ref{eq:cond2}).
{}From Proposition \ref{pr:main} and (\ref{eq:Le}), one has
\begin{theorem}\label{th:expand}
\begin{equation*}
L(u) = \sum_{i=n+2}^{N}(-1)^iT^{(N-i)}_1(u+\frac{i}{2})D^i - 
\sum_{i=0}^n(-1)^iT^{(i)}_1(u+\frac{i}{2})D^i.
\end{equation*}
\end{theorem}
{}From Proposition \ref{pr:Lproperty} and Theorem \ref{th:expand}
we conclude
$S_a\cdot T^{(b)}_1(u) = 0$ for all $1 \le a, b \le n$.

For later convenience, we extend $T^{(a)}_1(u)$ to $a \in \Z$ by
\begin{align}
&T^{(a)}_1(u) +T^{(N-a)}_1(u) = 0 \quad \forall a,\label{eq:T+T=0}\\
&T^{(a)}_1(u)=0 \quad \mbox{ for }\; a<0, \qquad T^{(0)}_1(u)=1.\nonumber
\end{align}
Then Theorem \ref{th:expand} is rephrased as
\begin{equation}\label{eq:LT}
L(u) = -\sum_{i=0}^{N}(-1)^iT^{(i)}_1(u+\frac{i}{2})D^i.
\end{equation}

Now we turn to the linear difference equation
\begin{equation}\label{eq:Lw=0}
L(u)w(u) = 0.
\end{equation}
Let $w_1(u), \ldots, w_N(u)$ be a basis of the solutions of (\ref{eq:Lw=0}).
For $i_1, \ldots, i_m \in \Z \; (m \le N)$, we prepare a 
shorthand for the Casorati determinant:
\begin{equation}\label{eq:[]def}
[i_1,\ldots,i_m] = \det\begin{pmatrix}
w_1(u+i_1) & \cdots & w_1(u+i_m) \\
\vdots & & \vdots \\
w_m(u+i_1) & \cdots & w_m(u+i_m)
\end{pmatrix}.
\end{equation}
We will write $[0,\ldots, 3]$ to mean the consecutive filling 
$[0,1,2,3]$ for example.
Note that $u$-dependence is suppressed in LHS and 
the overall shift $i_r \rightarrow i_r + 1$ is equivalent
to $u \rightarrow u+1$.
{}From (\ref{eq:Lw=0}) it follows that 
\begin{equation}\label{eq:shift}
[0,\ldots, N-1]  = -[1,\ldots, N].
\end{equation}
\begin{proposition}[Weyl type formula for fundamental $q$-characters]
\label{pr:wtype}
\[
T^{(a)}_1(u+\frac{a}{2}) = 
\frac{[0,\ldots, a-1, a+1,\ldots, N]}
{[1,2, \ldots, N]}\quad 0 \le a \le N.
\]
\end{proposition}
\begin{proof}
Solve the simultaneous equation obtained by setting 
$w = w_1,\ldots, w_N$ in (\ref{eq:Lw=0}), i.e., 
\begin{equation}\label{eq:L=0}
w(u+N) = \sum_{i=0}^{N-1}(-1)^iT^{(i)}_1(u+\frac{i}{2})w(u+i)
\end{equation}
with respect to the coefficients $T^{(i)}_1(u+\frac{i}{2})$.
\end{proof}

Proposition \ref{pr:wtype} is also shown
by applying Proposition \ref{pr:nnsy2} to Proposition \ref{pr:main}.
For the choice $w(u) = Q_1(u)$ (see (\ref{eq:LQ})), 
(\ref{eq:L=0}) is often 
called the ``$T-Q$ relation".

Before closing the section, we make 
a few miscellaneous remarks on $L(u)^{-1}$.
For $m \in {\mathbb Z}_{\ge 1}$ define 
\begin{equation*}
T^{(1)}_m(u+\frac{1}{2}) = \sum z_{i_1}(u-\frac{m-1}{2})
z_{i_2}(u-\frac{m-3}{2}) \cdots 
z_{i_m}(u+\frac{m-1}{2}),
\end{equation*}
where the sum runs over $i_1, \ldots, i_m \in J$ such that 
\begin{equation}\label{eq:cond3}
i_k \preceq i_{k+1} \; \mbox{ or }\;
(i_k,i_{k+1}) = (\bar{n},n)\; \mbox{ for } 1 \le k \le m-1.
\end{equation}
We put $T^{(0)}_1(u) = T^{(a)}_0(u) = 1$.
{}From the first expression in Lemma \ref{lem:kakikae}, 
we find
\begin{equation}\label{eq:Linv}
-L(u)^{-1} = \sum_{m=0}^\infty T^{(1)}_m(u+\frac{m}{2})D^m,
\end{equation}
hence $S_a \cdot T^{(1)}_m(u) = 0$.

The horizontal tableaux obeying the condition (\ref{eq:cond3}) 
have appeared in eq. (3.15a) of \cite{KS}.
We suppose that the quantity $T^{(1)}_m(u)$ 
is the irreducible $q$-character with highest weight monomial 
$\prod_{j=1}^mY_1\bigr(u+\frac{m+1-2j}{2}\bigl)$.
Multiplying (\ref{eq:LT}) and (\ref{eq:Linv}) we
deduce two types of ``$T-T$ relations":
\begin{align*}
&\sum_{a=0}^N(-1)^aT^{(1)}_{m-a}(u-\frac{a}{2})
T^{(a)}_1(u+\frac{m-a}{2}) = \delta_{m 0},
\\
&\sum_{a=0}^N(-1)^aT^{(1)}_{m-a}(u+\frac{m+a}{2})
T^{(a)}_1(u+\frac{a}{2}) = \delta_{m 0}
\end{align*}
for $m \in \Z_{\ge 0}$.
As is well known for $A_n$  case,  the $T-Q$ relation 
(\ref{eq:LQ}), namely, 
\begin{equation*}
\sum_{a=0}^N(-1)^aQ_1(u+a)T^{(a)}_1(u+\frac{a}{2}) = 0
\end{equation*}
is obtained {}from the limit 
$m \rightarrow \infty$ in either $T-T$ relations  
with a formal identification
\begin{equation*} 
Q_1(u) = \lim_{m \rightarrow \infty}T^{(1)}_m(u \mp \frac{m}{2}).
\end{equation*}
The $T-T$ relations are also obtainable by 
expanding the determinant expression for $T^{(1)}_m(u)$ 
in Remark \ref{rem:knh} with respect to the first row or $m$ th column.

\section{\mathversion{bold} 
Higher $q$-characters generated by $L(u)$}\label{sec:higher}
Evaluation of the ratio 
$\frac{[0,i_1, i_2, \cdots, i_{N-1}]}
{[0,1,2,\cdots, N-1]}$ 
of the Casorati determinants has been done in 
Appendix \ref{app:iran}.
Especially, Proposition \ref{pr:nnsy2} is an essential 
result saying that 
it is a polynomial in the fundamental $q$-characters 
$\{T^{(a)}_1(u) \mid 1 \le a \le n, u \in \C\}$.
To clarify its $q$-character content (irreducibility, 
decomposition into classical characters, etc.) 
for general $i_1, \ldots, i_{N-1}$ is left to a future study.
See Remark \ref{rem:pn}.
In this section we concentrate on a modest class generated 
by a repeated application of (\ref{eq:L=0}).
The resulting relation of the form 
\begin{equation}\label{eq:Dk}
w(u+k) = \sum_{i=0}^{N-1}(-1)^iH^{(i)}_{k}(u+\frac{i}{2})w(u+i)
\quad k \ge 0
\end{equation}
uniquely determines the coefficients.
In other words, we define 
$H^{(i)}_k(u)\; (0 \le i \le N-1, k \in \Z_{\ge 0}, u \in \C)$
via the recursion relation and the initial condition:
\begin{equation}\label{eq:recH}
\begin{split}
&H^{(i)}_{k+1}(u+\frac{i}{2}) = 
-T^{(i)}_1(u+\frac{i}{2})H^{(N-1)}_k(u+\frac{N+1}{2}) 
-H^{(i-1)}_k(u+\frac{i+1}{2}),\\
&H^{(i)}_0(u) = \delta_{i 0},
\end{split}
\end{equation}
where in the former we understand that $H^{(-1)}_k(u)=0$.
By the definition 
$H^{(i)}_k(u) = (-1)^i\delta_{i k}$ for $0 \le k \le N-1$ and 
$H^{(i)}_N(u) = T^{(i)}_1(u)$.
The following shows that 
the class $H^{(i)}_k(u)$ is relevant to Young diagrams of hook shape.
\begin{proposition}[Jacobi-Trudi and Weyl type formula]\label{pr:hook}
\begin{align*}
H^{(i)}_k(u+\frac{i}{2}) &= 
\begin{cases}\delta_{i k}(-1)^i & 0 \le k \le N-1\\
-\det_{1 \le j,l \le k-N+1}
\left(T^{(\la'_j-j+l)}_1\Bigl(u+\frac{N-2-\la'_j+j+l}{2}\Bigr)\right)
& k \ge N
\end{cases}\\
&=-\frac{[0,1,\ldots, i-1,i+1,i+2,\ldots,N-1,k]}{[0,\ldots,N-1]},
\end{align*}
where $\la'_j$ is specified {}from $i$ and $N$ by
$\la'_j = 1+(N-i-1)\delta_{j 1}$.
\end{proposition}
\begin{proof}
Due to (\ref{eq:T+T=0}), the determinant satisfies the same recursion as 
(\ref{eq:recH}), proving the first expression.
As for the second one, 
solve the simultaneous equation obtained by 
taking $w=w_1,\ldots, w_N$ in (\ref{eq:Dk}).
It can also be derived by rewriting the first one 
by using Proposition \ref{pr:nnsy2}.
\end{proof}
Regarding $\la'=(\la'_j)$ as the transpose of the Young diagram $\la$
(cf. \cite{M1}), we see that the latter is of hook shape 
with width $k-N+1$ and depth $N-i$.

Let us rewrite (\ref{eq:L=0}) and (\ref{eq:Dk}) into matrix forms.
We introduce the $N$-dimensional vectors and square matrices 
\begin{align*}
&\vec{w}(u) = 
\begin{pmatrix}
w(u) \\
w(u+1) \\
\vdots \\
w(u+N-1)
\end{pmatrix},
\qquad
\vec{h}_k(u) = 
\begin{pmatrix}
H^{(0)}_k(u) \\
-H^{(1)}_k(u+\frac{1}{2}) \\
\vdots \\
(-1)^{N-1}H^{(N-1)}_k(u+\frac{N-1}{2})
\end{pmatrix},\\
&{\mathcal T}(u) = \begin{pmatrix}
0 & 0 & \cdots & 0 & T^{(0)}_1(u) \\
1 & 0 &        & 0 & -T^{(1)}_1(u+\frac{1}{2}) \\
0 & 1 &        & \vdots & \vdots \\
\vdots & &     1  &  0 & (-1)^{N-2}T^{(N-2)}_1(u+\frac{N-2}{2}) \\
0 &  \cdots & 0 & 1 & (-1)^{N-1}T^{(N-1)}_1(u+\frac{N-1}{2})
\end{pmatrix},\\
&{\mathcal H}_k(u) = \left(
\vec{h}_k(u), \vec{h}_{k+1}(u), \ldots, \vec{h}_{k+N-1}(u) \right).
\end{align*}
Then (\ref{eq:L=0}) and (\ref{eq:Dk}) lead to
\begin{equation*}
\phantom{}^t \vec{w}(u+1) = \phantom{}^t \vec{w}(u){\mathcal T}(u), \qquad 
\phantom{}^t \vec{w}(u+k) = \phantom{}^t \vec{w}(u){\mathcal H}_k(u).
\end{equation*}
Therefore we have a product formula
\[
{\mathcal T}(u) {\mathcal T}(u+1) \cdots {\mathcal T}(u+k-1) = 
{\mathcal H}_k(u).
\]
See \cite{DDT,S1,S2} for a similar structure 
observed for the Stokes multipliers in $A^{(1)}_n$ case.
We set
\begin{equation}
\sigma_j = \begin{cases}
1 & 1 \le j \le n\\
-1 & n + 1 \le j \le N-1.
\end{cases} \label{eq:sigma}
\end{equation}

\begin{conjecture}\label{con:hook}
For $k \ge N+1$, the quantity 
$\sigma_i H^{(i)}_k(u+\frac{i}{2})\; (0 \le i \le N-1)$ is the 
irreducible $q$-character with the highest weight monomial
($Y_0(u) = 1$) 
\begin{equation*}
\begin{split}
&Y_n(u+\frac{n+2}{2}) \prod_{j=2}^{k-N}
Y_1\bigl(u+\frac{N+2j-1}{2}\bigr) \quad \mbox{ if } \;\; i = n+1,\\
&Y_{\min(i,N-i)}(u+\frac{i}{2}) \prod_{j=1}^{k-N}
Y_1\bigl(u+\frac{N+2j-1}{2}\bigr) \quad \mbox{ otherwise}.
\end{split}
\end{equation*}
\end{conjecture}

Let us turn to the $C_n$ content of $H^{(i)}_k$.
Let $(\alpha|\gamma)$  be the Young diagram of hook shape 
with width $\alpha+1$ and depth $\gamma+1$
(Frobenius notation \cite{M1}).
The corresponding $C_n$ character with highest 
weight $\alpha \La_1 + \La_{\gamma+1}$ is represented by
$\chi_{(\alpha|\gamma)}$.
Especially, we shall denote the character of the 
trivial representation by $\chi_{(-1|0)}=1$ 
rather than by $\chi_{(0|-1)}$.
Recall the homomorphism \cite{FR2}
\[
\beta: \; {\mathcal Y} = \Z[Y_a(u)^{\pm1}]_{1 \le a \le n, u \in \C} \rightarrow 
\Z[e^{\pm \La_a}]_{1 \le a \le n}
\]
sending $Y_a(u)^{\pm1}$ to $e^{\pm \La_a}$.
For $1 \le i \le N-1$ we know  
\[
\beta\bigl(\sigma_iT^{(i)}_1(u)\bigr)= \begin{cases}
0& \mbox{ if }\, i=n+1\\
\chi_{(0 \vert \min(i,N-i)-1)} &\mbox{ otherwise}.
\end{cases}
\]

\begin{proposition}\label{pr:hookchi}
For $k \ge N+1$, the image of 
$H^{(i)}_{k}(u)\; (0 \le i \le N-1)$ under $\beta$ is given as 
\begin{align*}
&\beta\bigl(H^{(0)}_k\bigr) = 
-\beta\bigl(H^{(N-1)}_{k-1}\bigr) = 
\sum^{\wedge} \chi_{(k-N-2j-1|0)},  \\
&\beta\bigl(H^{(a)}_k\bigr) = 
\begin{cases}
\sum^{\vee}  \chi_{(k-N-2j|a-1)}+
\sum^{\vee}  \chi_{(k-N-2j-1|a)} & 1 \le a \le n-1\\ 
-\sum^{\vee}  \chi_{(k-N-2j|N-a-1)}-
\sum^{\vee}  \chi_{(k-N-2j-1|N-a-2)} & n+2 \le a \le N-2,
\end{cases}\\
&\beta\bigl(H^{(n)}_{k-1}\bigr) = 
- \beta\bigl(H^{(n+1)}_k\bigr) =
\sum^{\vee}  \chi_{(k-N-2j-1|n-1)},
\end{align*}
where we have suppressed $u$ on LHS since the result is independent of it.
The sum  $\sum^{\wedge}\chi(\alpha-2j|\gamma)$ 
(resp. $\sum^{\vee}\chi(\alpha-2j|\gamma)$) 
extends over $j \in \Z_{\ge 0}$ such that 
$\alpha-2j \ge \min(0,\gamma-1)$
(resp. $\alpha-2j \ge 0$).
\end{proposition}
\begin{proof}
Check the relation (\ref{eq:recH}) under $\beta$ by means of 
($1 \le a \le n, \; p \ge 0$)
\begin{equation*}
\chi_{(p-1|0)}\chi_{(0|a-1)} = 
\chi_{(p|a-1)} + \chi_{(p-1|a)} +
\chi_{(p-1|a-2)} + \chi_{(p-2|a-1)},
\end{equation*}
where $\chi_{(p|a)}=0$ for any $p$ if $a\le -1$ or $a=n$.
\end{proof}

Extending $\beta$ on ${\mathcal Y}$, 
we introduce a map $\beta'$ which is applicable also to 
the solutions $w_i(u)$ of (\ref{eq:Lw=0}) as follows.
We parameterize the fundamental weight $\La_a$ 
in terms of the orthogonal basis 
$\veps_b\, (1 \le b \le n)$ as
$\La_a = \veps_1+ \cdots + \veps_a$.
The $C_n$ Weyl group acts as a 
permutation or $\times (\pm 1)$ on  $\{\veps_i\}$ as is well known.
For $1 \le i \le N$, introduce the variable $x_i$ by
\begin{equation*}
x_i = \begin{cases}
1 & i=n+1\\
-1 & i=n+2\\
e^{\sigma_i \veps_{\min(i,N+1-i)}} & \mbox{ otherwise}. 
\end{cases}
\end{equation*}
We take  $\beta'$ to be the same as $\beta$ on ${\mathcal Y}$, and 
\begin{equation*}
\beta':\; x_i(u) \mapsto x_i, \quad
w_i(u) \mapsto x^u_{N+1-i} \quad 1 \le i \le N.
\end{equation*}
The former ($i \neq n+1, n+2$) follows {}from 
(\ref{eq:defz})--(\ref{eq:xdef1}) by applying  $\beta$.
The latter has been adjusted to (\ref{eq:Ly=0}).
As a result we get
\begin{equation*}
\beta':\; \frac{[0,i_1, i_2, \cdots, i_{N-1}]}
     {[0,1,2,\cdots, N-1]}
\mapsto 
\frac{\det_{1 \le j,k \le N}\left(x^{i_{k-1}}_j\right)}
{\det_{1 \le j,k \le N}\left(x^{k-1}_j\right)},
\end{equation*}
where $i_0 = 0$.
This is a Weyl group invariant Laurent polynomial in 
$e^{\La_1}, \ldots, e^{\La_n}$, hence 
a linear combination of $C_n$ characters.
Under $\beta'$, Proposition \ref{pr:hook}
yields Jacobi-Trudi and Weyl type formulas 
(but $A_{2n+1}$-like rather than $C_n$) for the linear combinations 
of $C_n$ characters associated to the hook diagrams.

\begin{remark}\label{rem:pn}
For any $i_1, \ldots, i_{N-1} \in \Z$, one has 
$\frac{[0,i_1, i_2, \cdots, i_{N-1}]}{[0,1,2,\cdots, N-1]}
\in {\mathcal Y}$.
However the coefficients of the monomials 
in ${\mathcal Y}$ are not always all positive or negative. 
For example for $C_2$, one has 
\begin{align*}
\frac{[0,1,3,4,6,7]}{[0,1,2,3,4,5]}  &= 
T^{(1)}_1(u+\frac{3}{2})T^{(1)}_1(u+\frac{5}{2}) - 
T^{(2)}_1(u+1)T^{(2)}_1(u+3)\\
&= \frac{Y_1(u+\frac{7}{2})}{Y_1(u+\frac{11}{2})} 
-Y_2(u+2)Y_2(u+4) + Y_2(u+3) - \cdots,
\end{align*}
which consists of 19 monomials.
\end{remark}

\section{Solution of the $T$-system}\label{sec:Tsys}

The $T$-system is a set of functional relations 
among a certain family 
$\{T^{(a)}_m(u) \mid 1 \le a \le n, m \in \Z_{\ge 1}, u \in \C \}$
of commuting transfer matrices in solvable lattice models proposed 
in \cite{KNS} for any $U_q(X^{(1)}_n)$.
It has a form of the Toda field equation on a 
discrete space-time:
\[
T^{(a)}_m\left(u-\frac{(\alpha_a\vert \alpha_a)}{2}\right)
T^{(a)}_m\left(u+\frac{(\alpha_a\vert \alpha_a)}{2}\right) = 
T^{(a)}_{m+1}(u)T^{(a)}_{m-1}(u) + S^{(a)}_m(u),
\]
where $m \ge 1$, $T^{(a)}_0(u) = 1$ and 
$S^{(a)}_m(u)$ is a certain product of $T^{(a)}_m(u)$'s.
$\alpha_a\; (1 \le a \le n)$ denotes the simple root of $X_n$.
The $T$-system uniquely determines  $T^{(a)}_m(u)$ as 
a rational function of the 
fundamental ones $\{ T^{(a)}_1(u) \mid 1 \le a \le n, u \in \C \}$.
Moreover it has been proved for non-exceptional algebras \cite{KNH} that
$T^{(a)}_m(u)$ is actually a polynomial in these variables 
expressed as a determinant or pfaffian.

Let $W^{(a)}_m(u)$ denote the Kirillov-Reshetikhin module 
over $U_q(X^{(1)}_n)$ 
\cite{KR}.
By this we mean the irreducible one whose $q$-character 
has the highest weight monomial 
$\prod_{j=1}^m Y_a\bigl(u+\frac{(\alpha_a \vert \alpha_a)}{2}(m+1-2j)\bigr)$.
In the light of the correspondence between the transfer matrices and 
$q$-characters (cf. section 6.1 in \cite{FR2}), it is natural to make
\begin{conjecture}\label{con:Tsys}
The identification $T^{(a)}_m(u) = \chi_q(W^{(a)}_m(u))$ solves the $T$-system.
\end{conjecture}
Motivated by these aspects, we here present the solution 
of the $C^{(1)}_n$ $T$-system 
in terms of the Casorati determinants (\ref{eq:[]def}),
which may be viewed as a Weyl type formula for 
$q$-characters of the Kirillov-Reshetikhin module.
The $T$-system is explicitly given by
\begin{align}
T^{(a)}_m(u-\frac{1}{2}) T^{(a)}_m(u+\frac{1}{2})
 &= T^{(a)}_{m+1}(u) T^{(a)}_{m-1}(u) + T^{(a-1)}_{m}(u) T^{(a+1)}_{m}(u)
 \quad 1\le a\le n-2,  \label{t1} \\
T^{(n-1)}_{2m}(u-\frac{1}{2}) T^{(n-1)}_{2m}(u+\frac{1}{2})
 &= T^{(n-1)}_{2m+1}(u) T^{(n-1)}_{2m-1}(u) +
      T^{(n-2)}_{2m}(u)  T^{(n)}_{m}(u-\frac{1}{2})
T^{(n)}_{m}(u+\frac{1}{2}),
 \label{t2} \\     
T^{(n-1)}_{2m+1}(u-\frac{1}{2}) T^{(n-1)}_{2m+1}(u+\frac{1}{2})
 &= T^{(n-1)}_{2m+2}(u) T^{(n-1)}_{2m}(u) +
      T^{(n-2)}_{2m+1}(u)  T^{(n)}_{m}(u)  T^{(n)}_{m+1}(u),
 \label{t3} \\     
T^{(n)}_m(u-1) T^{(n)}_m(u+1) &= T^{(n)}_{m+1}(u) T^{(n)}_{m-1}(u) +
T^{(n-1)}_{2m}(u). \label{t4}
\end{align}
where $T^{(a)}_0(u)= T^{(0)}_m(u)=1$.

Set
\begin{align*}
\xi^{(a)}_m (u) &= [0,1,\cdots, a-1, a+m, a+m+1, \cdots, N+m-1],\\
\xi(u) &= \xi^{(1)}_0(u) = [0,1, \cdots, N-1].
\end{align*}

\begin{lemma}\label{lem:dual}
$\xi^{(a)}_m (u)$ satisfies the relations:
\begin{align}
&\xi^{(a)}_m(u)\xi^{(a)}_m(u+1) -
\xi^{(a)}_{m+1}(u)\xi^{(a)}_{m-1}(u+1) - \xi^{(a+1)}_{m}(u)\xi^{(a-1)}_{m}(u+1)=0,
\label{eq:plucker}\\
&\xi^{(a)}_m(u) = (-1)^{a-\frac{N}{2}+m} \xi^{(N-a)}_m (u+a -\frac{N}{2}).
\label{eq:dual}
\end{align}
\end{lemma}
\begin{proof}
(\ref{eq:plucker}) is a Pl\"ucker relation.
(\ref{eq:dual}) is shown by applying   
(\ref{eq:T+T=0}) to the latter formula in 
Proposition \ref{pr:nnsy2}.
\end{proof}
{}From $a=N/2$ case of (\ref{eq:dual}), it follows that
\begin{equation*}
\xi^{(n+1)}_m(u) = 0\quad m \in 2\Z_{\ge 0}+1. 
\end{equation*}
Consequently, one has the identities:
\begin{equation}\label{eq:comp}
\begin{split}
\xi^{(n)}_{2m}(u)  \xi^{(n+2)}_{2m}(u-1) &= 
\xi^{(n+1)}_{2m}(u)  \xi^{(n+1)}_{2m}(u-1),\\
\xi^{(n)}_{2m+1}(u)  \xi^{(n+2)}_{2m+1}(u-1) &= -
\xi^{(n+1)}_{2m}(u)  \xi^{(n+1)}_{2m+2}(u-1),
\end{split}
\end{equation}
since LHS $-$ RHS  become 
$-\xi^{(n+1)}_{2m+1}(u-1)  \xi^{(n+1)}_{2m-1}(u)$ and 
$\xi^{(n+1)}_{2m+1}(u)  \xi^{(n+1)}_{2m+1}(u-1)$ 
due to (\ref{eq:plucker}). 

Now our solution is given by
\begin{proposition}\label{pr:Tsys}
The following solves the $C^{(1)}_n$  $T$-system.
\begin{align}
&T^{(a)}_m(u+\frac{a+m-1}{2})=(-1)^m  \frac{\xi^{(a)}_m (u) }{\xi(u)}
\qquad 1 \le a \le n-1, \label{fort1} \\
&T^{(n)}_m(u+\frac{n+2m}{2}) T^{(n)}_m(u+\frac{n+2m-2}{2})
 = \frac{\xi^{(n)}_{2m} (u) }{\xi(u)},
  \label{fort2} \\
&T^{(n)}_m(u+\frac{n+2m}{2}) T^{(n)}_{m+1}(u+\frac{n+2m}{2})
 = \frac{\xi^{(n)}_{2m+1} (u) }{\xi(u+1)},
 \label{fort3}  \\
&T^{(n)}_m(u+\frac{n+2m}{2})^2 =
\frac{\xi^{(n+1)}_{2m} (u) }{\xi(u)}, 
 \label{fort4} 
\end{align}
which are equivalent, due to (\ref{eq:dual}), 
to the alternative forms:
\begin{align}
&T^{(a)}_m(u+\frac{a+m-1}{2})=
(-1)^{a-\frac{N}{2}} \frac{\xi^{(N-a)}_m (u+a-\frac{N}{2}) }{\xi(u)}
\quad 1 \le a \le n-1,  \label{fort1d} \\
&T^{(n)}_m(u+\frac{n+2m}{2}) T^{(n)}_m(u+\frac{n+2m-2}{2})
 =  \frac{\xi^{(n+2)}_{2m} (u-1) }{\xi(u+1)},  \label{fort2d} \\
&T^{(n)}_m(u+\frac{n+2m}{2}) T^{(n)}_{m+1}(u+\frac{n+2m}{2})
 =  \frac{\xi^{(n+2)}_{2m+1} (u-1) }{\xi(u+1)}.    \label{fort3d}
\end{align}
\end{proposition}
\begin{proof}
First we are to 
show the consistency of (\ref{fort2}), (\ref{fort2d}) and (\ref{fort4}).
Namely, evaluation of 
$(T^{(n)}_m(v) T^{(n)}_m(v-1) )^2$
by (\ref{fort2}) and (\ref{fort2d}) indeed coincides with
the result by (\ref{fort4}).
To see this, we multiply (\ref{fort2}) not with itself but with 
(\ref{fort2d}), leading to ($v=u+\frac{n+2m}{2}$)
\[
(T^{(n)}_m(v) T^{(n)}_m(v-1))^2 = 
\frac{\xi^{(n)}_{2m}(u)  \xi^{(n+2)}_{2m}(u-1) }{\xi(u) \xi(u-1)}.
\]
On the other hand (\ref{fort4}) says that LHS is equal to 
\[
\frac{\xi^{(n+1)}_{2m}(u)  \xi^{(n+1)}_{2m}(u-1) }{\xi(u) \xi(u-1)}.
\]
Thus they agree owing to the first relation in (\ref{eq:comp}).
Similarly, the consistency of (\ref{fort3}), (\ref{fort3d}) and 
(\ref{fort4}) is confirmed by means of the second relation in 
(\ref{eq:comp}).

Now we proceed to the check of (\ref{t1})--(\ref{t4}).
Upon substituting (\ref{fort1}), 
the difference of LHS and RHS of  (\ref{t1}) vanishes  
due to (\ref{eq:plucker}).
Similarly, one can verify (\ref{t2}) and (\ref{t3}) by using
(\ref{fort2}) and (\ref{fort3}).
As for the last relation (\ref{t4}), we first multiply 
$(T^{(n)}_m (u))^2 $ and regroup the factors as
\begin{equation*}
\begin{split}
&\left(T_m^{(n)}(u)T^{(n)}_m(u-1)\right)
\left(T^{(n)}_m(u+1)T^{(n)}_m(u)\right) \\
&= \left(T^{(n)}_m(u)T^{(n)}_{m+1}(u)\right)
\left(T^{(n)}_{m-1}(u)T^{(n)}_{m}(u)\right)+
\left(T^{(n)}_m(u)\right)^2 T^{(n-1)}_{2m}(u).
\end{split}
\end{equation*}
Upon applying (\ref{fort2}), (\ref{fort3}) and  (\ref{fort4}) to the first,
the second and the third terms, respectively, the result again 
reduces to the Pl\"ucker relation (\ref{eq:plucker}).
\end{proof}
Combining (\ref{fort2}), (\ref{fort3}), (\ref{fort2d}) and 
(\ref{fort3d}), we also have the expression
\[
T^{(n)}_m(u+\frac{n+2m}{2}) = (-1)^{m}\prod_{j=1}^m
\frac{\xi^{(n)}_{2j-1}(u+1)}{\xi^{(n)}_{2j-2}(u+1)}
= \prod_{j=1}^m
\frac{\xi^{(n+2)}_{2j-1}(u)}{\xi^{(n+2)}_{2j-2}(u)}.
\]
This implies the square root of (\ref{fort4}) is taken 
so that $T^{(n)}_m(u) = \prod_{j=1}^mY_n(u+m+1-2j) + \cdots$.

\begin{remark}\label{rem:knh}
The following Jacobi-Trudi type formula
is also known as Theorem 3.1 in \cite{KNH}: 
\begin{align*}
&T^{(a)}_m(u) = \det_{1 \le j,l \le m}\left(
T^{(a-j+l)}_1\bigl(u+\frac{j+l-m-1}{2}\bigr)\right)\quad 1 \le a \le n-1,\\
&T^{(n)}_m(u) = (-1)^m \mbox{pf}_{1 \le j,l \le 2m}
\left(
T^{(n+1-j+l)}_1\bigl(u+\frac{j+l-2m-1}{2}\bigr)\right),
\end{align*}
where $\mbox{pf}$ stands for the Pfaffian.
\end{remark}

\section{$B^{(1)}_n$ and $D^{(1)}_n$ cases}\label{sec:BD}

Here we present $L$ operators 
having the same property as Proposition \ref{pr:Lproperty} 
for $B^{(1)}_n$ and $D^{(1)}_n$.
However they are not polynomials in $D$.
Essentially they are the  
generating series of the pole-free combinations in the analytic Bethe ansatz
\cite{KOS,TK}.

For $B^{(1)}_n$ we set
\begin{align*}
&z_a(u) = \frac{Y_a(u+a)}{Y_{a-1}(u+a+1)}, \quad 
z_{\bar{a}}(u) = \frac{Y_{a-1}(u+2n-a)}{Y_{a}(u+2n-a+1)}
\quad 1 \le a \le n-1,\\
&z_n(u) = \frac{Y_n(u+\frac{2n+1}{2})Y_n(u+\frac{2n-1}{2})}
{Y_{n-1}(u+n+1)},\quad 
z_{\bar{n}}(u) = \frac{Y_{n-1}(u+n)}
{Y_n(u+\frac{2n+3}{2})Y_n(u+\frac{2n+1}{2})},\\
&z_0(u) = \frac{Y_n(u+\frac{2n-1}{2})}{Y_n(u+\frac{2n+3}{2})},
\end{align*}
{\small 
\begin{equation*}
L(u) = \displaystyle\stackrel{\longrightarrow}{\prod_{a=1}^n}
(1-z_{\bar{a}}(u)D^2)\cdot (1+z_{0}(u)D^2)^{-1}\cdot
\displaystyle\stackrel{\longleftarrow}{\prod_{a=1}^n}(1-z_a(u)D^2).
\end{equation*}
}

Let $\alpha_a \; (1 \le a \le n)$ be the simple root of $B_n$ 
normalized as $(\alpha_a \vert \alpha_a) = 2-\delta_{a n}$.
We let $S_a$ denote the screening operator for 
$B^{(1)}_n$ specified by 
(\ref{eq:SAS}) and (\ref{eq:defy}).

\begin{proposition}\label{pr:LB}
$S_a \cdot L(u) = 0$ for all $1 \le a \le n$.
\end{proposition}

\begin{proof}
For $a \neq n$ this can be directly checked as  (\ref{eq:SL}).
For $a=n$, we are to show $S_n \cdot X(u+n-2) = 0$, where
{\small 
\begin{equation*}
X(v) = 
\Big(1-\dfrac{Y_{n-1}(v+2)}{Y_{n}(v+5/2)Y_{n}(v+7/2)} D^2 \Big)
\Big(1+\dfrac{Y_{n}(v+3/2)}{Y_{n}(v+7/2)} D^2 \Big)^{-1}
\Big(1-\dfrac{Y_{n}(v+3/2)Y_{n}(v+5/2)}{Y_{n-1}(v+3)} D^2 \Big).
\end{equation*}
}
Expanding $X(v)$ one has
\begin{align*}
&X(v) = 1-f(v) D^2+h(v) \sum_{j=0}^{\infty}(-1)^j k(v+2 j) D^{2j+4},\\
&f(v) = \frac{Y_n(v+3/2)}{Y_n(v+7/2)} + 
\frac{Y_{n-1}(v+2)}{Y_n(v+5/2)Y_n(v+7/2)} +
\frac{Y_n(v+3/2)Y_n(v+5/2)}{Y_{n-1}(v+3)},\\
&k(v) = \frac{1}{Y_n(v+11/2)} + \frac{Y_n(v+9/2)}{Y_{n-1}(v+5)},\quad 
h(v) = Y_n(v+3/2)+ \frac{Y_{n-1}(v+2)}{Y_{n}(v+5/2)}.
\end{align*}
It is easy to verify 
$S_n\cdot f(v) = S_n\cdot k(v) 
= S_n\cdot h(v) = 0$.
\end{proof}
We introduce the expansion coefficients of $L(u)$ as
\begin{align}
L(u) &= 1+ \sum_{a \ge 1} (-1)^a T^a(u+a)D^{2a},\label{eq:ta}\\
L(u)^{-1} &= 1+ \sum_{m \ge 1}  T_m(u+m)D^{2m}.\label{eq:tm}
\end{align}
Under the correspondence (\ref{eq:defy}), 
they agree with those
defined in eq.(2.7) in \cite{KOS}.
For $1 \le a \le n-1$, $T^a(u)$ essentially coincides with 
$s_a(z)$ in section 11.2 of \cite{FR1}, which 
may be viewed as the $q$-character of the $a$ th 
fundamental representation of $U_q(B^{(1)}_n)$.
Note also that $L(u)Q_1(u) = 0$. 
As a corollary of Proposition \ref{pr:LB}, 
these coefficients are annihilated by all 
the screening operators $S_a$.

For $D^{(1)}_n$ we set 
\begin{align*}
&z_a(u) = \frac{Y_a(u+a)}{Y_{a-1}(u+a+1)}, \quad 
z_{\bar{a}}(u) = \frac{Y_{a-1}(u+2n-a-1)}{Y_{a}(u+2n-a)}
\quad 1 \le a \le n-2,\\
&z_{n-1}(u) = \frac{Y_n(u+n-1)Y_{n-1}(u+n-1)}
{Y_{n-2}(u+n)},\quad 
z_{ \,\overline{n-1}}\,(u) = \frac{Y_{n-2}(u+n)}{Y_{n}(u+n+1)Y_{n-1}(u+n+1)},\\
&z_n(u) = \frac{Y_n(u+n-1)}{Y_{n-1}(u+n+1)},\quad 
z_{\bar{n}}(u) = \frac{Y_{n-1}(u+n-1)}{Y_n(u+n+1)},
\end{align*}
{\small 
\begin{equation}\label{eq:LD}
L(u) = \displaystyle\stackrel{\longrightarrow}{\prod_{a=1}^n}
(1-z_{\bar{a}}(u)D^2)\cdot (1-z_{n}(u)z_{\bar{n}}(u+2)D^4)^{-1}\cdot
\displaystyle\stackrel{\longleftarrow}{\prod_{a=1}^n}(1-z_a(u)D^2).
\end{equation}
}

Let $\alpha_a \; (1 \le a \le n)$ be the simple root of $D_n$ 
normalized as $(\alpha_a \vert \alpha_a) = 2$.
We let $S_a$ denote the screening operator for 
$D^{(1)}_n$ specified by 
(\ref{eq:SAS}) and (\ref{eq:defy}).

\begin{proposition}\label{pr:LD}
$S_a \cdot L(u) = 0$ for all $1 \le a \le n$.
\end{proposition}
The proof is similar to Proposition \ref{pr:LB}.
In particular, $S_n\cdot L(u) =0$ is reduced to the 
following two lemmas.

\begin{lemma}\label{lem:w}
Setting
\begin{equation*}
h_a(u):= Y_a(u)+\frac{Y_{n-2}(u+1)}{Y_a(u+2)},\quad 
k_a(u):= Y_a(u)^{-1}+\frac{Y_a(u-2)}{Y_{n-2}(u-1)} \quad (a=n-1,n),
\end{equation*}
one has ${S}_n\cdot h_n(u) = {S}_n\cdot k_n(u) = 0$.
\end{lemma}

\begin{lemma}\label{lem:exp}
The $Y_n$-dependent factors in (\ref{eq:LD}) can be expanded as
{\small 
\begin{align*}
&\left(1-\frac{Y_{n-2}(v+4)}{Y_{n-1}(v+5)Y_n(v+5)}D^2\right)
\left(1-\frac{Y_{n-1}(v+3)}{Y_n(v+5)}D^2\right)
\left(1-\frac{Y_n(v+3)}{Y_n(v+7)}D^4\right)^{-1}\\
&\left(1-\frac{Y_n(v+3)}{Y_{n-1}(v+5)}D^2\right)
\left(1-\frac{Y_{n-1}(v+3)Y_n(v+3)}{Y_{n-2}(v+4)}D^2\right)\\
&=1-\sum_{j \ge 0}\left(k_{n-1}(v+4j+5)h_n(v+3) +
(1-\delta_{j,0})k_n(v+4j+5)h_{n-1}(v+3) \right)D^{4j+2}\\
&\quad +\sum_{j \ge 0}\left(k_{n-1}(v+4j+7)h_{n-1}(v+3) +
k_n(v+4j+7)h_n(v+3) -\delta_{j,0}
\frac{Y_{n-2}(v+4)}{Y_{n-2}(v+6)}\right)D^{4j+4},
\end{align*}
}
where $v= u+n-4$.
\end{lemma}
Defining $T^a(u)$ by  (\ref{eq:ta}) and (\ref{eq:LD}), one finds, 
under the correspondence (\ref{eq:defy}), that 
$T^a(u)$ coincides with ${\mathcal T}^a(u)$ in eq.(2.9) in \cite{TK}.
For $1 \le a \le n-2$, 
$T^a(u)$ essentially agrees with 
$s_a(z)$ in section 11.4 of \cite{FR1}, which 
may be viewed as the $q$-character of the $a$ th 
fundamental representation of $U_q(D^{(1)}_n)$.
Note also that $L(u)Q_1(u) = 0$. 
As a corollary of Proposition \ref{pr:LD}, 
$T^b(u)$ is annihilated by all 
the screening operators $S_a$.

\appendix
\section{Proof of Proposition \ref{pr:main}}\label{app:proof}

We begin by grouping the summands in RHS of (\ref{eq:tx}) as
\begin{equation*}
\sum_{i \not \ni n+1,n+2} + \sum_{i \ni n+1, i \not \ni n+2}
+ \sum_{i \ni n+2, i \not \ni n+1} + \sum_{i \ni n+1,n+2}.
\end{equation*}
Due to (\ref{eq:xdef2}), the 2nd and 3rd terms cancel each other.
The summands in the first and the last terms can be expressed 
with  $z_i(u)$'s by using (\ref{eq:xdef2}) and 
\begin{equation*}
x_{n+1}(u)x_{n+2}(u-1) = - z_n(u)z_{\bar{n}}(u-1).
\end{equation*}
The result reads
\begin{align}
&\sum_i X_{u+\frac{a}{2}-1}(i_1,\ldots,i_{a}) \nonumber \\
&= \sum_iZ_{u+\frac{a}{2}-1}(i_1,\ldots, i_a) - 
\sum_{k=0}^{a-2} \sum_{i,j}Z_{u+\frac{a}{2}-1}
(i_1,\ldots, i_k,n,\bar{n},j_1,\ldots, j_{a-2-k}).\label{eq:z-z}
\end{align}
Here the first sum is taken according to (\ref{eq:zsum}).
$\sum_{i,j}$ in the second term extends over 
$i_1,\ldots, i_k, j_1, \ldots, j_{a-2-k} \in J$
such that 
$1 \preceq i_1 \prec \cdots \prec i_k \preceq n$ and 
$\bar{n} \preceq j_1 \prec \cdots \prec j_{a-2-k} \preceq \bar{1}$.
Note that the second term contains the summands
with at most {\em two } $n$'s and {\em two} $\bar{n}$'s in the 
$i,j$-arrays.
So those patterns do not match a semistandard column tableau
with respect to the order $\prec$.
This will be the point 
that the remainder of the analysis concerns.
{}From (\ref{eq:z-z}) and (\ref{eq:cond1})--(\ref{eq:cond2}), 
Proposition \ref{pr:main} is reduced to 
($v = u + \frac{a}{2}-1$)
\begin{equation}\label{eq:reduced}
\begin{gathered}
\sum_{k=0}^{a-2} \sum_{i,j}Z_v
(i_1,\ldots, i_k,n,\bar{n},j_1,\ldots, j_{a-2-k})\\
= \sum_{1 \preceq i_1 \prec \cdots \prec i_a \preceq \bar{1}, \;
\text{ (\ref{eq:cond2}) is broken}}
Z_v(i_1,\ldots, i_a).
\end{gathered}
\end{equation}

For $n=2$ or $0 \le a \le 2$, 
it is straightforward to check (\ref{eq:reduced}).
Thus we assume $n \ge 3$ and fix
$3 \le a \le n$ (and any $v \in {\mathbb C}$) 
{}from now on.
We call the array $(i_1,\ldots, i_a)$ of elements 
$i_1,\ldots, i_a \in J$ a tableau.
We do not a priori assume (\ref{eq:cond1}) and (\ref{eq:cond2}).

\begin{lemma}\label{lem:tr}
\begin{equation*}
Z_v(\ldots,b,\overbrace{\ldots}^{n-b+1},\ol{b},\ldots) = 
Z_v(\ldots,b-1,\overbrace{\ldots}^{n-b+1},\ol{b-1},\ldots)
\quad 2 \le b \le n
\end{equation*}
with no change for $\ldots$ parts.
\end{lemma}
\begin{proof}
For any $u$ we have
\begin{equation}\label{eq:zeq}
z_b(u)z_{\bar{b}}(u-n+b-2) = z_{b-1}(u)z_{\ol{b-1}}(u-n+b-2) \quad
2 \le b \le n.
\end{equation}
\end{proof}
\noindent
Actually (\ref{eq:zeq}) is valid also for $b=1$ 
if  we interpret $z_0(u) = z_{\bar{0}}(u) = 1$.

We introduce a map of tableaux:
\begin{equation*}
\tau_b: \; (i_1, \ldots, i_a) \mapsto 
(i'_1,\ldots, i'_a) 
\quad 2 \le b \le n,
\end{equation*}
where RHS is obtained {}from LHS by making the transformation 
$(\ldots,b,\overbrace{\ldots}^{n-b+1},\ol{b},\ldots) \mapsto
(\ldots,b-1,\overbrace{\ldots}^{n-b+1},\ol{b-1},\ldots)$ 
for all the $(b,\ol{b})$ pairs matching this configuration.
When there is no such $(b,\bar{b})$ pair, 
we assume the action of $\tau_b$ is trivial, i.e.,
$(i'_,\ldots, i'_a) = (i_1,\ldots, i_a)$.
Due to  Lemma \ref{lem:tr}, $\tau_b$ is $Z_v$-preserving.
\begin{example}\label{ex:tau}
We set $n=a=9$. Omitting the parenthesis we have
\begin{align}
3\; 5 \; 7 \; 9 \; 9 \; \bar{9} \; \bar{8} \; \bar{7} \; \bar{3} 
\overset{\tau_9}\longmapsto 
&3\; 5 \; 7 \; 8 \; 9 \; \bar{8} \; \bar{8} \; \bar{7} \; \bar{3} \nonumber \\
\overset{\tau_8}\longmapsto 
&3\;  5 \; 7 \; 7 \; 9 \; \bar{8} \; \bar{7} \; \bar{7} \; \bar{3} \nonumber \\
\overset{\tau_7}\longmapsto 
&3\; 5 \; 6 \; 6 \; 9 \; \bar{8} \; \bar{6} \; \bar{6} \; \bar{3} \nonumber \\
\overset{\tau_6}\longmapsto 
&3\;  5 \; 5 \; 6 \; 9 \; \bar{8} \; \bar{6} \; \bar{5} \; \bar{3} \nonumber \\
\overset{\tau_5}\longmapsto 
&3\; 4 \; 5 \; 6 \; 9 \; \bar{8} \; \bar{6} \; \bar{4} \; \bar{3} 
\label{eq:koreya}\\
\overset{\tau_3}\longmapsto 
&2\;  4 \; 5 \; 6 \; 9 \; \bar{8} \; \bar{6} \; \bar{4} \; \bar{2}.
\label{eq:akan}
\end{align}
On all the  tableaux $\tau_2$ and $\tau_4$  act trivially.
\end{example}
Set 
\begin{align}
V &= \{(i_1\prec \cdots \prec i_k \preceq n,\bar{n} 
\preceq j_1 \prec \cdots \prec j_{a-2-k}) \mid 0 \le k \le a-2 \},
\label{eq:Vdef}\\
W &= \{ (i_1 \prec \cdots \prec i_a) \mid \text{(\ref{eq:cond2}) is broken} \},
\label{eq:Wdef}
\end{align}
where $i_r, j_r \in J$ extend over all the possibilities so that 
$V$ (resp. $W$) coincides with the range of the sum on  LHS (resp. RHS) 
of (\ref{eq:reduced}).
Since $\tau_b$'s are $Z_v$-preserving, Proposition \ref{pr:main} 
is reduced to constructing a bijection 
$\tau: V \longrightarrow W$ {}from their composition.
This will be achieved in Proposition \ref{pr:bij} later.

For $2 \le b \le n$ and $l,m \in {\mathbb Z}_{\ge 0}$ such that 
$\vert l-m \vert = 0, \pm1, (l,m) \neq (0,0)$, we introduce 
the subset $V^{l,m}_b \subset V$ defined by
\begin{align}
&V^{l,m}_b = \{
(i_1\prec \cdots \prec i_\alpha \prec 
\overbrace{b\cdots b}^l
\prec j_1 \prec \cdots \prec j_\beta \prec 
\overbrace{\ol{b}\cdots \ol{b}}^m
\prec k_1 \prec \cdots \prec k_\gamma) \mid 
\text{(A), (B), (C)}\},\label{eq:Vlmdef0}\\
&\text{(A) }\; i_r, j_r, k_r \in J,\; 
\alpha, \beta, \gamma \ge 0, \; \alpha + \beta + \gamma + l + m = a,
\label{eq:Vlmdef2}\\
&\text{(B) } \;  \text{if } j_r = d, j_s = \ol{d} 
\text{ for some  } b < d \le n, \text{ then } n+r-s \ge d, 
\label{eq:Vlmdef1}\\
&\text{(C)  One of the following (C1)-(C4) holds:} 
\label{eq:Vlmdef3}\\
&\quad \text{(C1): } l=m \ge 1, \; l+\beta = n-b+1, \text{ i.e., } \; 
(\ldots,\overbrace{\underline{b\ldots b}\; b}^l,\ldots,
\overbrace{\ol{b}\;\underline{\ol{b}\ldots\ol{b}}}^l,\ldots), \nonumber\\
&\quad \text{(C2): } l=m \ge 1, \; l+\beta = n-b+2, \text{ i.e., } \; 
(\ldots,\overbrace{\underline{b\ldots b}}^l,\ldots,
\overbrace{\underline{\ol{b}\ldots\ol{b}}}^l,\ldots),\nonumber\\
&\quad \text{(C3): } l=m+1 \ge 1, \; l+\beta = n-b+2, \text{ i.e., } \; 
(\ldots,\overbrace{\underline{b\ldots b}\; b}^l,\ldots,
\overbrace{\underline{\ol{b}\ldots\ol{b}}}^{l-1},\ldots),\nonumber\\
&\quad \text{(C4): } l=m-1 \ge 0, \; l+\beta = n-b+1, \text{ i.e., } \; 
(\ldots,\overbrace{\underline{b\ldots b}}^l,\ldots,
\overbrace{\ol{b}\; \underline{\ol{b}\ldots\ol{b}}}^{l+1},\ldots),\nonumber
\end{align}
where the underlines indicate those $b$ and $\ol{b}$ that are
changed into $b-1$ and $\ol{b-1}$ under the action of $\tau_b$ for 
$b \ge 2$.
In short, the condition (C) is selecting the tableaux of the form
\begin{equation}\label{eq:short}
(\ldots,\underline{b\ldots b}(b),\ldots,
(\ol{b}) \underline{\ol{b}\ldots \ol{b}},\ldots),
\end{equation}
where $(b)$ and $(\ol{b})$ can be present or absent 
independently.
The condition (B) says that (\ref{eq:cond2}) is satisfied for the 
segment $j_1 \prec \cdots \prec j_\beta$.
In (\ref{eq:Vlmdef0}), the inequalities involving 
$b$ or $\ol{b}$ should be imposed even when $l=0$ or $m=0$.
We set  
\begin{equation*}
V_b = \bigsqcup_{l,m \ge 0,\; (l,m) \neq (0,0)} V^{l,m}_b
\qquad 2 \le b \le n,
\end{equation*}
where we assume $V^{l,m}_b = \emptyset$ 
unless  $\vert l - m \vert = 0,\pm 1$. 

As an example, we have
\begin{equation*}
V = V_n = V^{1,1}_n \sqcup V^{1,2}_n \sqcup V^{2,1}_n \sqcup V^{2,2}_n,
\qquad V^{1,1}_n \subset W. 
\end{equation*}
Given $3 \le a \le n$, only 
$V_n, V_{n-1}, \ldots, V_{n-a+2}$ are non-empty, 
and the last one reads 
\begin{align*}
&V_{a^\sharp} = V^{1,1}_{a^\sharp} \sqcup V^{0,1}_{a^\sharp} 
\sqcup V^{1,0}_{a^\sharp}, \qquad 
{a^\sharp} = n-a+2,  \\
&V^{1,1}_{a^\sharp} = \{({a^\sharp} 
\prec j_1 \prec \cdots \prec j_{a-2} \prec \ol{{a^\sharp}})\},  \\
&V^{0,1}_{a^\sharp} = \{(j_1 \prec \cdots \prec j_{a-1} 
\prec \ol{{a^\sharp}})\} \quad (a^\sharp \prec j_1), \\
&V^{1,0}_{a^\sharp} = \{({a^\sharp} \prec j_1 \prec \cdots \prec j_{a-1})\}
\quad (j_{a-1}\prec  \ol{{a^\sharp}}),
\end{align*}
where $j_r$'s obey the condition (B).
It is easy to see
\begin{equation}\label{eq:end}
\tau_{a^\sharp}(t) = t \quad\text{for any }\; t \in V_{a^\sharp}.
\end{equation}
Note that $a^\sharp \ge 2$.
\begin{lemma}\label{lem:step}
If $t \in V_b$, then $\tau_b(t) = t$ or $\tau_b(t) \in V_{b-1}$ for 
$3 \le b \le n$.
\end{lemma}
\begin{proof} We set $c=b-1$.
As in (\ref{eq:short}), an element $t \in V_b$ has the form
$t = (\ldots,d,\underline{b\ldots b}(b),\ldots,
(\ol{b}) \underline{\ol{b}\ldots \ol{b}},e,\ldots)$,
where $d \prec b, \ol{b} \prec e$.
Suppose $\tau_b(t) \neq t$. 
We classify the nontrivial action of $\tau_b$ into four cases; 
(D1) $d=c, e=\ol{c}$, (D2) $d=c, e \succ \ol{c}$,
(D3) $d \prec c, e=\ol{c}$ and (D4) $d \prec c, e \succ \ol{c}$.
In each case, the action $t \rightarrow \tau_b(t)$ is given as follows:
\begin{alignat*}{2}
&\text{(D1): } (\ldots,c,\underline{b\ldots b}(b),\ldots,
(\ol{b})\underline{\ol{b}\ldots\ol{b}}\; \ol{c},\ldots)\;\; &\mapsto
(\ldots,\underline{c\; c\ldots c}(b),\ldots,
(\ol{b})\underline{\ol{c}\ldots\ol{c}\;\ol{c}},\ldots), \\
&\text{(D2): } (\ldots,c,\underline{b\ldots b}(b),\ldots,
(\ol{b})\underline{\ol{b}\ldots\ol{b}},\ldots) &\mapsto
(\ldots,\underline{c\; c\ldots c}\;c (b),\ldots,
(\ol{b})\underline{\ol{c}\ldots\ol{c}},\ldots), \\
&\text{(D3): } (\ldots,\underline{b\ldots b}(b),\ldots,
(\ol{b})\underline{\ol{b}\ldots\ol{b}}\; \ol{c},\ldots) &\mapsto
(\ldots,\underline{c\ldots c}(b),\ldots,
(\ol{b})\ol{c} \;\underline{\ol{c}\ldots\ol{c} \; \ol{c}},\ldots), \\
&\text{(D4): } (\ldots,\underline{b\ldots b}(b),\ldots,
(\ol{b})\underline{\ol{b}\ldots\ol{b}},\ldots) &\mapsto
(\ldots,\underline{c\ldots c}\; c (b),\ldots,
(\ol{b})\ol{c}\;\underline{\ol{c}\ldots\ol{c}},\ldots). 
\end{alignat*}
Here the underlines on RHS (resp. LHS) 
designate those $c, \ol{c}$ (resp. $b, \ol{b}$) 
changed under $\tau_c$ (resp. $\tau_b$).
Comparing (\ref{eq:short}) with RHS's of (D1)-(D4), 
we see that $\tau_b(t)$ satisfies the condition (C) (\ref{eq:Vlmdef3}) 
for $V_c$.
The condition (A) is clear.
The condition (B) is nontrivial only for the pair 
$(b,\ol{b})$ when $(b)$ and $(\ol{b})$ are
both present.
In such cases 
the number $\beta$ of the letters between $(b)$ and 
$(\ol{b})$ in (D1)-(D4) is bounded by 
$\beta \le n-b-1$.
Thus $n-\beta \ge b+1 > b$, showing that (B) is valid.
\end{proof}

Given $t=t_{n} \in V = V_n$, let $t_d = \tau_{d+1}\cdots \tau_n(t)$.
{}From (\ref{eq:end}) and Lemma \ref{lem:step}, we deduce 
\begin{lemma}\label{lem:pexist}
There exists a unique $p$ such that 
\begin{align*}
&a^\sharp \le p \le n, \\
&t_{p} \neq t_{p+1} \neq \cdots \neq t_{n},\quad 
 t_d \in V_d\; \text{ for }\;  p \le d \le n,\\
&t_p = t_{p-1}.
\end{align*}
\end{lemma}
\begin{lemma}\label{lem:inW}
$t_p \in W$.
\end{lemma}
\begin{proof}
In (D1)-(D4), the length of the underlines on the RHS are shorter 
than those on LHS only for (D4).
Therefore the situation 
\begin{equation*}
t_{p+1} \overset{\tau_{p+1}}{\longmapsto} t_p
\overset{\tau_{p}}{\longmapsto} t_{p-1} = t_p
\end{equation*}
means that the map $\tau_{p+1}$ underwent the pattern (D4) 
with length one underlines in its LHS. 
Namely we have 
\begin{equation}\label{eq:tttnamely}
\begin{gathered}
(\ldots, d,\; \underline{p+1}\; 
\overbrace{(p+1),\ldots,(\ol{p+1})}^{n-p}\; 
\underline{\ol{p+1}},\; e,\ldots) = t_{p+1} \in V_{p+1}\\
\overset{\tau_{p+1}}{\longmapsto}
(\ldots, d, \; p \; 
\overbrace{(p+1),\ldots,(\ol{p+1})}^{n-p}\; 
\ol{p},\; e,\ldots) = t_p \in V_p,
\end{gathered}
\end{equation}
where $d \preceq p-1$ and $\overline{p-1} \preceq e$.
Let us check that RHS belongs to $W$.
(Irrespectively  
of the presence or absence of $(p+1)$ and $(\ol{p+1})$, 
(\ref{eq:tttnamely}) says that 
there are always $n-p$ letters between $p$ and $\ol{p}$.)
First, there is no repetition of the same letter 
in the tableau $t_p$ because of the definition of $V_p$.
Second, the $(p,\ol{p})$ pair in the center certainly 
breaks the condition (\ref{eq:cond2}).
\end{proof}
\begin{remark}\label{rem:inside}
By the definition (\ref{eq:Wdef}), any tableau in $W$ contains 
a pair $(q,\ol{q})$ breaking (\ref{eq:cond2}).
Let us call such a pair with the largest value of $1 \le q \le n$ 
the {\em maximal breaking } pair.
(Actually a pair $(1,\ol{1})$ can never break (\ref{eq:cond2}).)
In the end of the proof of Lemma \ref{lem:inW}, 
we have also established the following:
the maximal breaking pair of $t_p$ is $(p,\ol{p})$.
\end{remark}
\begin{lemma}\label{lem:mbp}
Suppose that in the tableau
\begin{equation*}
(\ldots,q,\overbrace{\ldots}^{\gamma},\ol{q},\ldots) \in W,
\end{equation*}
the pair $(q,\ol{q})$ is the maximal breaking one.
Then we have $\gamma = n-q$.
\end{lemma}
\begin{proof}
Since $(q,\ol{q})$ breaks (\ref{eq:cond2}), we know $q \ge n-\gamma$.
Assume that $q > n-\gamma$, hence
$\sharp\{q+1,q+2,\ldots,n\} = n-q < \gamma$.
Then there is at least one pair $(r,\ol{r})$ with $q < r \le n$ between 
$q$ and $\ol{q}$, hence the tableau looks as
\begin{equation*}
(\ldots,q,\overbrace{\ldots}^{\alpha},r,
\overbrace{\ldots}^{\delta},\ol{r},
\overbrace{\ldots}^{\beta}\ol{q},\ldots), \quad 
\gamma = \alpha + \beta + \delta + 2.
\end{equation*}
By the definition, the pair $(r,\ol{r})$ must satisfy
the condition (\ref{eq:cond2}), meaning $r < n-\delta$.
When there are more than one such $r$, 
we take the smallest one among those, which implies 
$\alpha + \beta \le r-q-1$.
Now these relations lead to the contradiction:
\[
0< q-n+\gamma = q-n+\alpha+\beta+\delta+2 \le
-n+\delta+r+1 \le 0.
\]
\end{proof}
Note the consistency of (\ref{eq:tttnamely}), Remark \ref{rem:inside}
and Lemma \ref{lem:mbp}.

By virtue of 
Lemma \ref{lem:pexist} and Lemma \ref{lem:inW}, we are 
entitled to define
\begin{equation*}
\begin{split}
\tau : V &\longrightarrow W \\
t &\longmapsto t_p = \tau_{p+1}\tau_{p+2}\cdots \tau_n(t),
\end{split}
\end{equation*}
where $p$ is specified in Lemma \ref{lem:pexist}.
In Example \ref{ex:tau}, LHS is an element of $V$.
When calculating its image under $\tau$, one has $p=4$, and the 
answer is (\ref{eq:koreya}) and not (\ref{eq:akan}).
Observe that  $(4,\ol{4})$ is certainly the 
maximal breaking pair in (\ref{eq:koreya}) 
containing  $n-4=5$ letters in between.
\begin{proposition}\label{pr:bij}
The map $\tau: \; V \longrightarrow W$ is a bijection.
\end{proposition}
\begin{proof}
We have only to construct the inverse of $\tau$.
For $3 \le b \le n$ we define the map 
$\sigma_b: (i_1,\ldots, i_a) \mapsto (i'_1,\ldots,i'_a)$ 
by making the transformation 
$(\ldots,b-1,\overbrace{\ldots}^{n-b+1},\ol{b-1},\ldots)
\mapsto 
(\ldots,b,\overbrace{\ldots}^{n-b+1},\ol{b},\ldots)$ 
for all the $(b-1,\ol{b-1})$ pairs matching this configuration.
Given any tableau $s \in W$, 
we set $\sigma(s) = \sigma_n \cdots \sigma_{p+2}\sigma_{p+1}(s)$, where
$p$ is determined {}from the condition that 
$(p,\ol{p})$ is the maximal breaking pair of $s$.
By construction it is then evident that 
$\sigma(\tau(t))=t,\; \tau(\sigma(s)) = s$ for any 
$t \in V$ and $s \in W$.
\end{proof}

\section{Basic lemmas}\label{app:iran}

We keep the notation (\ref{eq:[]def}) but 
do not assume (\ref{eq:Lw=0}) and (\ref{eq:shift}) 
in Lemma \ref{lem:nnsy}.
Set
\begin{equation}\label{eq:xtdef}
\tilde{x}_m(u) = \frac{[0,\ldots,m-1][2,\ldots,m]}
{[1,\ldots,m][1,\ldots,m-1]}\quad 1 \le m \le N.
\end{equation}
We define $\tilde{e}_a(u)$ by (\ref{eq:edef}) by replacing   
$x_m(u)$ with $\tilde{x}_m(u)$.
In particular, $\tilde{e}_a(u) = 0$ if $a > N$ or $a<0$.
\begin{lemma}[\cite{NNSY}]\label{lem:nnsy}
Given the integers  $0=i_0 <i_1< \cdots < i_{N-1}$, 
let $\mu= (\mu_j)$ be the Young diagram with 
depth less than $N$ specified by $\mu_j = i_{N-j}+j - N$.
Let $\mu'=(\mu'_j)$ denote the transpose of $\mu$.
Assume $m \ge i_{N-1}-N+1$. Then
\begin{align*}
\frac{[0,i_1, i_2, \cdots, i_{N-1}]}
     {[m,\ldots, m+N-1]}
&= \sum_{t}\prod_{(\alpha,\beta) \in (m^N)/\mu}{\tilde x}_{t(\alpha,\beta)}
(u+\alpha+\beta-2)\\
&= \det_{1 \le j,l \le m}\left(
\tilde{e}_{N-\mu'_j-l+j}\Bigl(u+\frac{N-2+j+l-\mu'_j}{2}\Bigr)\right),
\end{align*}
where the sum $\sum_t$ extends over the semistandard tableaux 
on the skew Young diagram $(m^N)/\mu$ \cite{M1} on letters $\{1,\ldots, N\}$.
$t(\alpha,\beta)$ denotes the entry of $t$ at the 
$\alpha$th row and the $\beta$th column {}from the bottom left corner.
\end{lemma}
The lemma is related to the ninth variation of the Schur 
function \cite{M2}, and applicable to 
$q$-characters for $U_q(A^{(1)}_{N-1})$.
It is actually valid for any $N \in \Z_{\ge 1}$.
To approach the $U_q(C^{(1)}_n)$ case in question, 
we next take the constraints (\ref{eq:Lw=0}) 
and (\ref{eq:shift}) into account.
We introduce the difference operators $L_j(u)\; (1 \le j \le N)$  by
\begin{equation}\label{eq:defLj}
L_j(u) = \displaystyle\stackrel{\longrightarrow}{\prod_{i=N+1-j}^{N}}
(D-\epsilon_ix_i(u+n+1-i)),
\end{equation}
where $\epsilon_i=1$ except  $\epsilon_{n+1}=\epsilon_{n+2}=-1$.
By Lemma \ref{lem:kakikae} and (\ref{eq:Ldef}) 
we have $L(u) = L_N(u)$.
Take the basis $\{ w_1(u),\ldots, w_N(u) \}$
of the solutions to (\ref{eq:Lw=0}) such that 
\begin{equation}\label{eq:Ly=0}
L_j(u)w_m(u) = 0 \quad 
1 \le m \le j \le N.
\end{equation}
\begin{lemma}\label{lem:x=xt}
Under the above choice of the basis, we have
${\tilde x}_m(u) = x_m(u)$, where the latter is defined in 
(\ref{eq:xdef1})--(\ref{eq:xdef2}).
\end{lemma}

\begin{proof}
By calculating (\ref{eq:defLj}) directly, one gets $(1 \le j \le N-1)$ 
\begin{align}
L_j(u) &= D^j + (-1)^j\sigma'_j \frac{q_j(u+1)}{q_j(u)} 
+ \text{ terms proportional to } D,\ldots, D^{j-1}, 
\label{eq:const}\\
\sigma'_j &= \begin{cases}
1 & 1 \le j \le n+1\\
-1 & n + 2 \le j \le N-1,
\end{cases} \nonumber \\
q_j(u) &= \begin{cases}
Q_j(u+\frac{j-1}{2}) & 1 \le j \le n-1\\
Q_n(u+\frac{n}{2})Q_n(u+\frac{n-2}{2}) & j = n\\
Q_n(u+\frac{n}{2})^2 & j=n+1\\
Q_n(u+\frac{n}{2})Q_n(u+\frac{n+2}{2}) & j = n+2\\
Q_{N-j}(u+\frac{j-1}{2}) & n+3 \le j \le N-1.
\end{cases}\nonumber
\end{align}
Let $[0,\ldots,j-1]$ be the Casorati determinant of 
$w_m(u)$'s as defined in (\ref{eq:[]def}).
Due to (\ref{eq:Ly=0}) and (\ref{eq:const}) it satisfies 
the first order linear difference equation
\begin{equation*}
\left(D - \sigma'_j \frac{q_j(u+1)}{q_j(u)}\right)
[0,\ldots,j-1] = 0\quad 1 \le j \le N-1.
\end{equation*}
Thus we may set 
\begin{equation}\label{eq:cq}
[0,\ldots,j-1] = \phi_j(u)q_j(u),
\end{equation}
where $\phi_j(u)$ is any function satisfying
$\phi_j(u+1) = \sigma'_j\phi_j(u)$.
Substituting (\ref{eq:cq}) into (\ref{eq:xtdef}) one finds
${\tilde x}_m(u) = x_m(u)$. 
\end{proof}

Due to Lemma \ref{lem:x=xt}, we may set $\tilde{e}_a(u) = e_a(u)$.
By combining Proposition \ref{pr:main}, (\ref{eq:esym}) and (\ref{eq:T+T=0}),
this can be further identified with $T^{(a)}_1(u)$ for all $a \in \Z$.
Substituting this back to Lemma \ref{lem:nnsy} and 
using (\ref{eq:T+T=0}), we obtain
\begin{proposition}\label{pr:nnsy2}
Let $i_0,\ldots, i_{N-1}, \mu$ and $\mu'$  be 
as in Lemma \ref{lem:nnsy}.
Assume further that $w_1, \ldots, w_N$ satisfy (\ref{eq:Ly=0}).
Then we have 
\begin{align*}
\frac{[0,i_1, i_2, \cdots, i_{N-1}]}
     {[0,\ldots, N-1]}
&= (-1)^{\mu_1}
\sum_{t}\prod_{(\alpha,\beta) \in (\mu_1^N)/\mu}x_{t(\alpha,\beta)}
(u+\alpha+\beta-2)\\
&= \det_{1 \le j,l \le \mu_1}\left(
T^{(\mu'_j-j+l)}_1\Bigl(u+\frac{N-2+j+l-\mu'_j}{2}\Bigr)\right),
\end{align*}
where the sum $\sum_t$ extends over the semistandard tableaux 
on the skew Young diagram $(\mu_1^N)/\mu$.
$t(\alpha,\beta)$ is the entry of $t$ at the 
$\alpha$th row and the $\beta$th column {}from the bottom left corner
of the skew Young diagram $(\mu_1^N)/\mu$.
\end{proposition}

\medskip

\section*{Acknowledgements}
A.K., M.O. and Y.Y. have been supported by Grant-in-Aid for 
Scientific Research {}from the Ministry of Education, Culture, Sports,
Science and Technology of Japan.
The work of J.S. has been supported by
a Grand-in-Aid for Encouragement of Young Scientists
from the Japanese Society for the Promotion of Science, $\sharp$12740244.


\begin{thebibliography}{A}

\bibitem [B]{B}
R.J.~Baxter,
\textit{Exactly solved models in statistical mechanics}, Academic Press, 
London (1982).

\bibitem [BHK]{BHK}
V. V. Bazhanov, A. N. Hibberd and  S. M. Khoroshkin,
\textit{Integrable structure of $W_3$ conformal field theory, 
quantum Boussinesq theory and boundary affine toda theory},
hep-th/0105177.

\bibitem [CK]{CK} V. Chari and M. Kleber,
\textit{Symmetric functions and representations of quantum affine algebras},
math.QA/0011161.

\bibitem [CP1]{CP1} V.\ Chari and A.\ Pressley,
\textit{ A guide to quantum groups},
Cambridge Univ. Press. Cambridge (1994).


\bibitem [CP2]{CP2} V.\ Chari and A.\ Pressley,
\textit{ Quantum affine algebras and their
representations},
Canadian Math.\ Soc.\ Conf.\ Proc.\ {\bf 16}
(1995) 59--78.

\bibitem [DDT]{DDT} P. Dorey, C. Dunning and R. Tateo,
\textit{Differential equations for general SU(n) Bethe ansatz system},
J. Phys. A. Math. Gen. {\bf 33} (2000) 8427-8442.


\bibitem [FM]{FM} E. Frenkel and E. Mukhin,
\textit{Combinatorics of $q$-characters of finite dimensional 
representations of quantum affine algebras},
math.QA/9911112.

\bibitem [FR1]{FR1} E. Frenkel and N. Reshetikhin,
\textit{Quantum affine algebras and deformations of the Virasoro and 
${\mathcal W}$-algebras}, Comm. Math. Phys. {\bf 178} (1996) 237--264. 


\bibitem [FR2]{FR2} E. Frenkel and N. Reshetikhin,
\textit{The $q$-characters of representations of quantum affine 
algebras and deformations of ${\mathcal W}$-algebras},
Contemporary Math.\ {\bf 248} (1999) 163--205.


\bibitem [FRS]{FRS} E. Frenkel, N. Reshetikhin and M. A. Semenov-Tian-Shansky,
\textit{Drinfeld-Sokolov reduction for difference operators
and deformations of ${\mathcal W}$-algebras. I. The case of Virasoro algebra},
Comm. Math. Phys. {\bf 192} (1998)  605--629. 


\bibitem [KR]{KR}
A.\ N.\ Kirillov and N.\ Yu.\ Reshetikhin,
\textit{Representations of Yangians and multiplicity of
occurrence of the irreducible components of the
tensor product of representations of simple Lie algebras},
J.\ Sov.\ Math.\ {\bf 52} (1990) 3156--3164.


\bibitem[Ka]{Ka}
M.~Kashiwara, 
\textit{On level zero representations of quantized affine algebras},
math.QA/0010293.


\bibitem [Kn]{Kn} H. Knight,
\textit{Spectra of tensor products of finite-dimensional representations 
of Yangians},  J. Algebra {\bf 174} (1995) 187--196. 


\bibitem [KLWZ]{KLWZ} 
I. Krichever, O. Lipan, P. Wiegmann and A. Zabrodin, 
\textit{Quantum integrable models and discrete classical Hirota equations},
Comm. Math. Phys. {\bf 188} (1997) 267--304. 


\bibitem [KNH]{KNH} A. Kuniba, S. Nakamura and R. Hirota, 
\textit{Pfaffian and determinant 
solutions to a Discretized Toda equation 
for $B_r, C_r$ and $D_r$},
J. Phys. A: Math. Gen. {\bf 29} (1996) 1759--1766.


\bibitem [KNS]{KNS} A. Kuniba, T. Nakanishi and J. Suzuki,
\textit{Functional relations in solvable lattice models:
I. Functional relations and representation theory},
Int. J. Mod. Phys. A {\bf 9} (1994) 5215--5266.


\bibitem [KOS]{KOS} A. Kuniba, Y. Ohta and J. Suzuki,
\textit{Quantum Jacobi-Trudi and Giambelli formulae for 
$U_q(B^{(1)}_r)$ {}from the analytic Bethe ansatz},
J. Phys. A. Math. Gen. {\bf 28} (1995) 6211--6226.


\bibitem [KS]{KS}  A. Kuniba and J. Suzuki,
\textit{Analytic Bethe ansatz for fundamental representations of Yangians},
Comm. Math. Phys. \ {\bf 173} (1995) 225--264.

\bibitem [M1]{M1} I. G. Macdonald,
\textit{Symmetric functions and Hall polynomials}, 2nd ed.
Oxford Univ. Press, Oxford (1995).

\bibitem [M2]{M2} I. G. Macdonald,
\textit{Schur functions: Theme and variations},
Publ. I.R.M.A. Strasbourg, Acte $28^e$,
S\`eminaire Lotharingien, (1992) 5--39.

\bibitem [NNSY]{NNSY} J. Nakagawa, M. Noumi, M. Shirakawa and Y. Yamada,
\textit{Tableau representation for Macdonald's ninth variation of 
Schur functions}, preprint.


\bibitem [R]{R} N. Yu. Reshetikhin,
\textit{The spectrum of the transfer matrices connected with 
Kac-Moody algebras}, Lett. Math. Phys. {\bf 14} (1987) 235--246.

\bibitem [SS]{SS} M. A. Semenov-Tian-Shansky and A. V. Sevostyanov,
\textit{Drinfeld-Sokolov reduction for difference operators and
deformations of ${\mathcal W}$-algebras. II. The general semisimple case},
Comm. Math. Phys. {\bf 192} (1998) 631--647.

\bibitem [S1]{S1} J. Suzuki,
\textit{Functional relations in Stokes multipliers
and solvable models related to $U_q(A^{(1)}_n)$},
J. Phys. A. Math. Gen.{\bf 33} (2000) 3507-3521.


\bibitem [S2]{S2} J. Suzuki,
\textit{Stokes Multipliers, Spectral Determinants and $T-Q$ relations},
preprint.


\bibitem [TK]{TK} Z. Tsuboi and A. Kuniba,
\textit{Solutions of a discretized Toda field equation for $D_r$
{}from analytic Bethe ansatz},
J. Phys. A. Math. Gen. {\bf 29} (1996) 7785--7796.
\end{thebibliography}
\end{document}